\documentclass[preprint]{elsarticle}

\usepackage{lineno,hyperref}
\usepackage{amsmath, amssymb,epsfig}
\usepackage{enumerate} 
\usepackage{multirow} 
\usepackage{placeins}
\modulolinenumbers[5]

\journal{Journal of Computational and Applied Mathematics}

\newtheorem{thm}{Theorem}

\newdefinition{rmk}{Remark}
\newdefinition{dfn}{Definition}
\newdefinition{ex}{Example}
\newproof{pf}{Proof}

\newcommand{\ds}{\displaystyle} 

\counterwithin*{equation}{section}

\begin{document}

\begin{frontmatter}

\title{Modified Trapezoidal Product Cubature Rules.
Definiteness, Monotonicity and a Posteriori Error
Estimates\tnoteref{This research was supported by the Bulgarian
National Research Fund under Contract KP-06-N62/4.}}
\tnotetext[mytitlenote]{This research was supported by the Bulgarian
National Research Fund under Contract KP-06-N62/4.}



\author{Geno Nikolov\corref{mycorrespondingauthor}}
\cortext[mycorrespondingauthor]{Corresponding author}
\ead{geno@fmi.uni-sofia.bg}

\author{Petar Nikolov}
\ead{petartv@yahoo.com}

\address{Faculty of Mathematics and Informatics,
         Sofia University St Kliment Ohridski,
         5~James Bourchier Blvd., 1164 Sofia, Bulgaria}

\begin{abstract}
We study two modifications of the trapezoidal product cubature
formulae, approximating double integrals over the square domain
$[a,b]^2=[a,b]\times [a,b]$. Our modified cubature formulae use
mixed type data: except evaluations of the integrand on the points
forming a uniform grid on $[a,b]^2$, they involve two or four
univariate integrals. An useful property of these cubature formulae
is that they are definite of order $(2,2)$, that is, they provide
one-sided approximation to the double integral for real-valued
integrands from the class
$$
\mathcal{C}^{2,2}[a,b]=\{f(x,y)\,:\,\frac{\partial^4 f}{\partial
x^2\partial y^2}\ \text{continuous and does not change sign in}\
(a,b)^2\}.
$$
For integrands from $\mathcal{C}^{2,2}[a,b]$ we prove monotonicity
of the remainders and derive a-posteriori error estimates.
\end{abstract}

\begin{keyword}
Peano kernel representation\sep trapezium quadrature rule\sep
product cubature formulae\sep interpolation by blending
functions\sep a posteriori error estimates

\MSC[2010] 41A55\sep 65D30\sep 65D32
\end{keyword}

\end{frontmatter}


\section{Introduction and statement of the main result}

A standard tool for approximation of the definite integral $\ds
\int\limits_a^b f(x)\,dx$ is the ($n$-th) composite trapezium rule
$$
Q_n^{Tr}[f]=h_n\sideset{}{'}\sum\limits_{i=0}^{n}f(x_{i,n}), \quad
x_{i,n}=a+ih_n,\ h_n=\frac{b-a}{n}
$$
(here and in what follows, the notation $\sideset{}{'}\sum\nolimits$
means that the first and the last summands are halved). Denote by
$R[Q_n^{Tr};f]$ the remainder functional,
$$
R[Q_n^{Tr};f]:=\int\limits_a^b f(x)\,dx-Q_n^{Tr}[f]\,.
$$
Assuming the integrand $f$ is convex or concave in $[a,b]$, we have
the following well-known properties of the trapezium rule:
\begin{itemize}
\item[(i)] \emph{Definiteness:} $Q_n^{Tr}[f]$ provides one-sided
approximation to $I[f]$, namely,
$$
R[Q_n^{Tr};f]\leq 0\quad\text{for convex }f\;,\quad
R[Q_n^{Tr};f]\geq 0\quad\text{for concave }f\;;
$$
\item[(ii)] \emph{Monotonicity:}
$$
\big|R[Q_{2n}^{Tr};f]\big|\leq
\frac{1}{2}\big|R[Q_{n}^{Tr};f]\big|\,;
$$
\item[(iii)] \emph{A posteriori error estimate:}
$$
\big|R[Q_{2n}^{Tr};f]\big|\leq\big|Q_n^{Tr}[f]-Q_{2n}^{Tr}[f]\big|\,.
$$
\end{itemize}
Definiteness and monotonicity are properties of quadrature rules
which have been intensively studied (see e.g., \cite{AN2017,
KJF:1993a, KJF:1993, FKN:1998, KJ:1976, KN:1995, DN:1974, GN:1992,
GN:1996, GS:1972} and the monographs \cite{HB:1977} and
\cite{BP:2011}). These properties have found application in the
software packages for numerical integration.

The product trapezium cubature rule
$$
C_n^{Tr}[f]=h_n^2\sideset{}{'}\sum\limits_{i=0}^{n}\sideset{}{'}
\sum\limits_{j=0}^{n}f(x_{i,n},y_{j,n}), \qquad
x_{i,n}=y_{i,n}=a+ih_n,\ \ h_n=\frac{b-a}{n},
$$
is a natural candidate for approximating the double integral
$$
I[f]:=\iint\limits_{[a,b]^2}f(x,y)\,dxdy\,,
$$
where $[a,b]^2$ is the square $[a,b]\times [a,b]\subset
\mathbb{R}^2$. Since in the univariate case the convexity/concavity
of $f$ is secured by the assumption $f\in C^2[a,b]$ with
$f^{\prime\prime}\geq 0$ (resp. $f^{\prime\prime}\leq 0$) in
$(a,b)$, an appropriate class of functions in the bivariate case is
$$
\mathcal{C}^{2,2}[a,b]=\{f\in C^{2,2}([a,b]^2)\,:\,D^{2,2}f\
\text{does not change sign in}\ (a,b)^2\}\,.
$$
Here and henceforth, for $k,l\in\mathbb{N}$,
$$
C^{k,l}([a,b]^2):=\{f\,:\,[a,b]^2\mapsto \mathbb{R}, D^{i,j}f\
\text{is cont. in } [a,b]^2,\ 0\leq i\leq k,\ 0\leq j\leq l\},
$$
where
$$
D^{i,j}f:=\frac{\partial^{i+j} f}{\partial x^i\partial y^j}.
$$

In general, the remainder functional
$$
R\big[C_n^{Tr};f\big]=I[f]-C_n^{Tr}[f]
$$
does not possess properties like (i), (ii) and (iii) for integrands
$f\in\mathcal{C}^{2,2}[a,b]$. By suitably modifying $C_n^{Tr}$, we
obtain two families of cubature formulae which admit properties
analogous to (i), (ii) and (iii). We pay a price  for this: our
modified cubature formulae involve evaluations of two or four
univariate integrals of some traces of the integrand. With every
bivariate function $f(x,y)$ defined on $[a,b]^2$ we associate six
univariate functions:
\begin{eqnarray*}
&&f_l(t)=f(a,t),\quad f_r(t)=f(b,t),\quad f_d(t)=f(t,a),\quad
f_u(t)=f(t,b),\\ &&f_v(t)=f\big(\frac{a+b}{2},t\big),\quad
f_h(t)=f\big(t,\frac{a+b}{2}\big).
\end{eqnarray*}
\begin{dfn}
For every $n\in \mathbb{N}$, the modified product trapezoidal rules
$S_n^{-}$ and $S_n^{+}$ are defined as follows:
\begin{equation}\label{e:1.1}
S_n^{-}[f]=C_n^{Tr}[f]+(b-a)\Big(R[Q_n^{Tr};f_v]+R[Q_n^{Tr};f_h]\Big),
\end{equation}
\begin{equation}\label{e:1.2}
S_n^{+}[f]=C_n^{Tr}[f]\!+\!
\frac{b-a}{2}\Big(R[Q_n^{Tr};f_l]\!+\!R[Q_n^{Tr};f_r]\!+\!R[Q_n^{Tr};f_d]
\!+\!R[Q_n^{Tr};f_u]\Big).
\end{equation}
\end{dfn}
The cubature formulae $S_n^{+}$ and $S_n^{-}$ have been introduced
(with $[a,b]=[-1,1]$ and somewhat different notation) in
\cite{GN2001}, where some properties of these cubature formulae have
been established. Set
$$
R\big[S_n^{-};f\big]=I[f]-S_n^{-}[f], \quad
R\big[S_n^{+};f\big]=I[f]-S_n^{+}[f].
$$
The following result has been proved in \cite{GN2001}  (see Theorems
4.4 and 5.2 therein):

\smallskip\noindent
\textbf{Theorem A.} \emph{
\begin{enumerate}[~~(i)~~]
\item
For every $f\in C^{2,2}([a,b]^2)$, there exist $P\in [a,b]^2$ such
that
\begin{equation}\label{e:1.3}
R[S_n^{-};f]=c(S_n^-)\,D^{2,2}f(P),\quad
c(S_n^-)=-\frac{(b\!-\!a)^{6}}{144n^2}\Big(1\!+\!\frac{1}{n^2}\Big).
\end{equation}
\item
For every $f\in C^{2,2}([a,b]^2)$, there exist $P\in [a,b]^2$ such
that
\begin{equation}\label{e:1.6}
R[S_n^{+};f]=c(S_n^+)\,D^{2,2}f(P),\quad
c(S_n^+)=\frac{(b-a)^{6}}{72n^2}\Big(1-\frac{1}{2n^2}\Big).
\end{equation}
\end{enumerate}}
\smallskip

Our main result here is monotonicity and a posteriori error
estimates for the remainders of $S_n^{-}$ and $S_n^{+}$.
\begin{thm}\label{t:1.1}
If $f\in\mathcal{C}^{2,2}[a,b]$, then
\begin{equation}\label{e:1.4}
\big|R[S_{2n}^{-};f]\big|\leq \frac{1}{2}\,\big|R[S_{n}^{-};f]\big|
\end{equation}
and
\begin{equation}\label{e:1.5}
\big|R[S_{2n}^{-};f]\big|\leq
\big|S_{2n}^{-}[f]-S_{n}^{-}[f]\big|\,.
\end{equation}
\end{thm}
\begin{thm}\label{t:1.2}
If $f\in\mathcal{C}^{2,2}[a,b]$, then
\begin{equation}\label{e:1.7}
\big|R[S_{2n}^{+};f]\big|\leq
\Big(\frac{1}{2}+\frac{1}{4(2n-1)}\Big) \big|R[S_{n}^{+};f]\big|
\end{equation}
and
\begin{equation}\label{e:1.8}
\big|R[S_{2n}^{+};f]\big|\leq \frac{4n-1}{4n-3}\,
\big|S_{2n}^{+}[f]-S_{n}^{+}[f]\big|\,.
\end{equation}
\end{thm}

Let us briefly comment on the results stated in the above theorems.
Set
\begin{eqnarray*}
&&\mathcal{C}^{2,2}_{+}[a,b]=\{f\in \mathcal{C}^{2,2}[a,b]\,:\,
D^{2,2}f\geq 0\ \text{ in }\ [a,b]^2\},\\
&&\mathcal{C}^{2,2}_{-}[a,b]=\{f\in \mathcal{C}^{2,2}\,:\, -f\in
\mathcal{C}^{2,2}_{+}[a,b]\}.
\end{eqnarray*}
Theorem~A implies
$$
S_n^{+}[f]\leq I[f] \leq S_{m}^{-}[f], \qquad
f\in\mathcal{C}^{2,2}_{+}[a,b]
$$
and the reversed inequalities when $f\in\mathcal{C}^{2,2}_{-}[a,b]$.
Thus, for integrands from $\mathcal{C}^{2,2}[a,b]$ we obtain
inclusions for the true value of $I[f]$. By analogy with the
univariate case, we say that $S_n^{-}$ and $S_n^{+}$ are
respectively \emph{negative definite} and \emph{positive definite}
cubature formulae of order $(2,2)$.

Inequalities \eqref{e:1.4} and \eqref{e:1.7} in Theorems~\ref{t:1.1}
and \ref{t:1.2} imply that for integrands from
$\mathcal{C}^{2,2}[a,b]$ doubling $n$ results in reducing the error
magnitude by a factor of at least two (or almost two). Inequalities
\eqref{e:1.5} and \eqref{e:1.8} provide a posteriori error
estimates, which may serve as stoping rules (rules for termination
of calculations) in routines for automated numerical integration.
Let us point out here that, although $S_n^{\pm}$ involve univariate
integrals, the error bounds in \eqref{e:1.5} and \eqref{e:1.8} are
in terms of only point evaluations of the integrand. In a certain
sense, the inequalities in Theorem~\ref{t:1.1} and
Theorem~\ref{t:1.2} are the best possible, as it will become clear
from their proof.

The rest of the paper is structured as follows. In section 2 we
provide some necessary background, including Peano kernel
representation of linear functionals, in particular, of the
remainders of quadrature formulae, and the approach for construction
of cubature formulae using mixed type data through interpolation by
blending functions. At the end of this section we prove
Theorem~\ref{t:4}, which provides a method for derivation of a
posteriori error estimates for definite cubature formulae. In
section 3 we apply Theorem~\ref{t:4} to prove Theorem~\ref{t:1.1}
and Theorem~\ref{t:1.2}. In section~4 we present a numerical example
illustrating our theoretical results.
\section{Preliminaries}
\subsection{Peano kernel representation of linear functionals}
The notation $\pi_{m}$ will stand for the set of algebraic
polynomials of degree not exceeding $m$. A classical result of Peano
\cite{GP:1913} asserts that if $\mathcal{L}$ is a linear functional
defined in $C[a,b]$ which vanishes on $\pi_{r-1}$, then for
functions $g$ with $g^{(r-1)}$ absolutely continuous in $[a,b]$,
$\mathcal{L}[g]$ admits the integral representation
$$
\mathcal{L}[g]=\int\limits_{a}^{b}K_r(t)g^{(r)}(t)\,dt,
$$
where
$$
K_r(t)=\mathcal{L}\Big[\frac{(\cdot-t)_{+}^{r-1}}{(r-1)!}\Big],\quad
u_{+}=\max\{u,0\}.
$$
By $W^{r}_p[a,b]$, $\,(r\in \mathbb{N},\ p\geq 1)$, we denote the
Sobolev class of functions
$$
W^{r}_p[a,b]:=\{f\in C^{r-1}[a,b]\,:\, f^{(r-1)} \mbox{ abs.
continuous},\; \int_{a}^{b}\!|f^{(r)}(t)|^p\,dt<\infty\}\,.
$$
Note that $C^{r}[a,b]\subset W^{r}_p[a,b]$ for every $p\geq 1$.

A quadrature formula $Q_n$ approximating the integral $\ds
\int\limits_a^b f(x)\,dx$,
\begin{equation*}
Q_n[f]=\sum_{i=1}^{n}a_{i,n}\,f(\tau_{i,n})\,,\quad a\le
\tau_{1,n}<\tau_{2,n}<\cdots<\tau_{n,n}\leq b,
\end{equation*}
is said to have \emph{algebraic degree of precision} $r-1$ (in
short, $ADP(Q_n)=r-1$), if its remainder
$$
R[Q_n;f]=\int\limits_{a}^{b}f(x)\,dx-Q_n[f]
$$
vanishes identically on $\pi_{r-1}$ and $R[Q_n;x^{r}]\ne 0$. By
Peano's representation theorem, if $ADP(Q_n)=r-1$ and the integrand
$f$ belongs to $W^{s}_1[a,b]$, where $s\leq r$, $s\in\mathbb{N}$,
then the remainder of $Q_n$ admits the integral representation
\begin{equation}\label{e:2.1}
R[Q_n;f]= \int_{a}^{b}K_s(Q_n;t)f^{(s)}(t)\,dt
\end{equation}
with
$$
K_s(Q_n;t)=R\Big[Q_n;\frac{(\cdot-t)_{+}^{s-1}}{(s-1)!}\Big].
$$
$K_s[Q_n;\cdot]$ will be referred to as \emph{the $s$-th Peano
kernel of} $Q_n$. Explicit representations of $\,K_s(Q_n;t)\,$ for
$t\in [a,b]$ are
\begin{equation}\label{e:2.2}
\begin{split}
K_s(Q_n;t)&=(-1)^{s}\Big[\frac{(t-a)^s}{s!}-
\frac{1}{(s-1)!}\sum_{i=1}^{n}a_{i,n}(t-\tau_{i,n})_{+}^{s-1}\Big]\\
&=\frac{(b-t)^s}{s!}-
\frac{1}{(s-1)!}\sum_{i=1}^{n}a_{i,n}(\tau_{i,n}-t)_{+}^{s-1}\,.
\end{split}
\end{equation}
Application of H\"{o}lder inequality to \eqref{e:2.1} yields sharp
error estimates of the form
$$
\big|R[Q_n;g]\big|\leq c_{s,p}(Q_n)\,\Vert g^{(s)}\Vert_p,\quad
s=1,\ldots, r,
$$
in the corresponding Sobolev classes of functions, where
$\Vert\cdot\Vert_p$ is the $L_p([a,b])$-norm, $1\leq p\leq\infty$.

Of particular interest are the so-called definite quadrature
formulae. Quadrature formula $Q_n$ is said to be definite of order
$r$, if there exists a real non-zero constant $c_{r}(Q_n)$ such that
for the remainder of $Q_n$ there holds
$$
R[Q_n;f]=c_r(Q_n)\,f^{(r)}(\xi)
$$
for every $f\in C^{r}[a,b]$, with some $\xi\in [a,b]$ depending on
$f$. Furthermore, $Q_n$ is called positive definite (resp., negative
definite) of order $r$, if $c_r(Q_n)>0$ ($c_r(Q_n)<0$). It is clear
from \eqref{e:2.1} that $Q_n$ is a positive (negative) definite
quadrature formula of order $r$ if and only if $ADP(Q_n)=r-1$ and
$K_r(Q_n;t)\geq 0$ (resp., $K_r(Q_n;t)\leq 0$) on $(a,b)$. The
interest in definite quadratures of order $r$ lies in the one-sided
approximation they provide to the integral for integrands having
$r$-th derivative with a permanent sign in $[a,b]$.

We refer to \cite{HB:1977} and \cite[Chapt.4]{BP:2011} for more
details on the Peano kernel theory.
\subsection{Cubature formulae using mixed type of data}
In \cite{GN2001} V. Gushev and one of the authors initiated study on
adopting Peano kernel theory to the error estimation of cubature
formulae for approximate evaluation of double integrals on a
rectangular region, a further development can be found in
\cite{GN2009}. The approach applied there makes use of bivariate
interpolation by blending functions (see, e.g., \cite{BDHN:1999}).
We present below the main idea and some of the results from
\cite{GN2001} and \cite{GN2009}.

For $m_1,\,m_2\in\mathbb{N}$, the set of blending functions
$B^{m_1,m_2}([a,b]^2)$ is defined by
\begin{equation}\label{e2.3}
B^{m_1,m_2}([a,b]^2)=\big\{f\in C^{m_1,m_2}([a,b]^2)\;:\;
D^{m_1,m_2}f=0\big\}.
\end{equation}
Given $X=\{x_1,x_2,\ldots,x_{m_1}\}$, $Y=\{y_1,y_2,\ldots,y_{m_2}\}$
with $a\!\leq\! x_1<\cdots<x_{m_1}\!\leq\! b$ and $a\leq
y_1<\cdots<y_{m_2}\leq b$, we define a blending grid $G=G(X,Y)$ by
$$
G(X,Y)=\Big\{(x,y)\;:\;\
\prod_{\mu=1}^{m_1}(x-x_\mu)\prod_{\nu=1}^{m_2}(y-y_\nu)=0\Big\}.
$$
For any function $f$ defined on $[a,b]^2$ there exists a unique
(Lagrange) blending interpolant $Bf=B_Gf\in B^{m_1,m_2}([a,b]^2)$
satisfying $Bf_{|G(X,Y)}=f_{|G(X,Y)}$, it is explicitly given by
$$
Bf=\mathcal{L}_x f+\mathcal{L}_y f-\mathcal{L}_x \mathcal{L}_y f,
$$
where $\mathcal{L}_x$ and $\mathcal{L}_y$ are the Lagrange
interpolation operators with respect to $x$ and $y$ with sets of
interpolation nodes $X$ and $Y$, respectively. Precisely, if
$\{\ell_{\mu}\}_{\mu=1}^{m_1}$ and
$\{\widetilde{\ell}_{\nu}\}_{\nu=1}^{m_2}$ are the Lagrange basic
polynomials for $\pi_{m_1-1}$ and $\pi_{m_2-1}$, defined by
$\ell_{\mu}(x_j)=\delta_{\mu,j}\ (j=1,\ldots,m_1)$ and
$\widetilde{\ell}_{\nu}(y_j)=\delta_{\nu,j}\ (j=1,\ldots,m_2)$,
respectively, with $\delta_{i,j}$ being the Kronecker symbol, then
$$
Bf(x,y)=\sum_{\mu=1}^{m_1}\ell_{\mu}(x)f(x_\mu,y)+
\sum_{\nu=1}^{m_2}\widetilde{\ell}_{\nu}(y)f(x,y_\nu)-
\sum_{\mu=1}^{m_1}\sum_{\nu=1}^{m_2}
\ell_{\mu}(x)\widetilde{\ell}_{\nu}(y)f(x_{\mu},y_{\nu}).
$$
\begin{rmk}
Blending grid $G=G(X,Y)$ can be defined also by sets $X$ and $Y$
containing multiple nodes, then $\mathcal{L}_x$ and $\mathcal{L}_y$
are replaced by the associated Hermite interpolation operators.
\end{rmk}
For $f\in C^{r,s}([a,b]^2)$ with $r,\,s\in \mathbb{N}$ satisfying
$r\leq m_1$ and $s\leq m_2$, two iterated applications of the Peano
theorem to $f-Bf=(Id-\mathcal{L}_x)(Id-\mathcal{L}_y)f$, where $Id$
is the identity operator, imply
\begin{equation}\label{e:2.4}
f(x,y)-Bf(x,y)= \iint\limits_{[a,b]^2} \mathcal{K}_r(x,t)
\widetilde{\mathcal{K}}_s(y,\tau)D^{r,s}f(t,\tau)\,dtd\tau,
\end{equation}
where
$$
\mathcal{K}_r(x,t)=\frac{1}{(r-1)!}\,\Big((x-t)_{+}^{r-1}
-\sum_{\mu=1}^{m_1}\ell_{\mu}(x)(x_{\mu}-t)_{+}^{r-1}\Big),
$$
$$
\widetilde{\mathcal{K}}_s(y,\tau)=\frac{1}{(s-1)!}\,\Big((y-\tau)_{+}^{s-1}
-\sum_{\nu=1}^{m_2}\widetilde{\ell}_{\nu}(y)(y_{\nu}-\tau)_{+}^{s-1}\Big).
$$
Integrating \eqref{e:2.4} over $[a,b]^2$ we obtain the approximation
$$
I[f]\approx I[Bf]=:C[f].
$$
$C[f]$ is called \emph{blending cubature formula} and is of the form
\begin{equation}\label{e:2.5}
C[f]=\sum_{\mu=1}^{m_1}b_{\mu}\int\limits_{a}^{b}f(x_{\mu},t)\,dt +
\sum_{\nu=1}^{m_2}\widetilde{b}_{\nu}\int\limits_{a}^{b}f(t,y_{\nu})\,dt
-\sum_{\mu=1}^{m_1}\sum_{\nu=1}^{m_2}b_{\mu}\widetilde{b}_{\nu}f(x_{\mu},y_{\nu}),
\end{equation}
where $Q^{\prime}[g]=\sum_{\mu=1}^{m_1}b_{\mu}g(x_{\mu})$ and
$Q^{\prime\prime}[g]=\sum_{\nu=1}^{m_2}\widetilde{b}_{\nu}g(x_{\nu})$
are the interpolatory quadrature formulae with sets of nodes $X$ and
$Y$, respectively. The remainder of this blending cubature formula
admits the integral representation
\begin{equation}\label{e:2.6}
R[C;f]=I[f]-C[f]=\iint\limits_{[a,b]^2}K_r(Q^{\prime},t)K_s(Q^{\prime\prime};\tau)
D^{r,s}f(t,\tau)\,dtd\tau.
\end{equation}
Now, assuming that the integrand belongs to a certain Sobolev class
of bivariate functions, one can derive error estimates using
H\"{o}lder's inequality.

The presence of $m_1+m_2$ univariate integrals in the blending
cubature formula \eqref{e:2.5} is a drawback, as the precise values
of these integrals may be not accessible. A sequence of blending
cubature formulae involving only these $m_1+m_2$ univariate
integrals and approximating $I[f]$ with increasing accuracy is
constructed through the scheme
$$
I[f-Bf]\approx C_n[f-Bf],
$$
where $\{C_n\}$ is a sequence of cubature formulae of standard type,
i.e., using only point evaluations. We denote the resulting cubature
formulae by $S_n$,
\begin{equation}\label{e:2.7}
S_n[f]=I[Bf]+C_n[f]-C_n[Bf].
\end{equation}
A reasonable choice for $C_n$ is to be a product cubature formula
generated by a quadrature formula
\begin{equation}\label{e:2.8}
Q_n[g]=\sum_{j=0}^{n}a_{j,n}g(x_{j,n})\approx\int\limits_{a}^{b}g(t)\,dt,
\end{equation}
i.e.
$$
C_n[f]=\sum_{i=0}^{n}\sum_{j=0}^{n}a_{i,n}a_{j,n}f(x_{i,n},x_{j,n}).
$$
With such a choice of $C_n$, \eqref{e:2.4} and \eqref{e:2.5} imply
an integral representation of the remainder of $S_n$,
\begin{equation}\label{e:2.9}
R[S_n,f]=I[f]-S_n[f]=
\iint\limits_{[a,b]^2}K_{r,s}(S_n;t,\tau)D^{r,s}f(t,\tau)\,dtd\tau
\end{equation}
with
$$
K_{r,s}(S_n;t,\tau)=K_r(Q^{\prime};t)K_s(Q^{\prime\prime};\tau)
-Q_n[\mathcal{K}_r(\cdot,t)]Q_n[\widetilde{\mathcal{K}}_{s}(\cdot,\tau)].
$$
Further representations of $S_n$ and its Peano kernel
$K_{r,s}(S_n;\cdot)$ are given in the following theorem:
\begin{thm}[\cite{GN2001}, Theorem 3.1]\label{t:3}
Assume that in the blending cubature formula $S_n$ defined by
\eqref{e:2.7}, $C_n$ is the product cubature formula generated by a
quadrature formula $Q_n$ with $ADP(Q_n)\geq\max\{ADP(Q^{\prime}),
ADP(Q^{\prime\prime})\}$, where $Q^{\prime}$ and $Q^{\prime\prime}$
are the interpolatory quadrature formulae with sets of nodes $X$ and
$Y$, respectively. Then for $(t,\tau)\in [a,b]^2$ the Peano kernel
$K_{r,s}(S_n;t,\tau)$ has the following representations:
\begin{equation}\label{e:2.10}
K_{r,s}(S_n;t,\tau)=K_r(Q^{\prime};t)K_s(Q_n;\tau)
+K_r(Q_n;t)Q_n[\widetilde{\mathcal{K}}_s(\cdot,\tau)],
\end{equation}
\begin{equation}\label{e:2.11}
K_{r,s}(S_n;t,\tau)=K_s(Q^{\prime\prime};\tau)K_r(Q_n;t)
+K_s(Q_n;\tau)Q_n[\mathcal{K}_r(\cdot,t)],
\end{equation}
\begin{equation}\label{e:2.12}
\begin{split}
K_{r,s}(S_n;t,\tau)=&K_r(Q^{\prime};t)K_s(Q_n;\tau)
+K_s(Q^{\prime\prime};\tau))K_r(Q_n;t)\\
&-K_r(Q_n;t))K_s(Q_n;\tau).
\end{split}
\end{equation}
The cubature formula $S_n$ admits the representation
\begin{equation}\label{e:2.13}
S_n[f]=C_n[f]+Q^{\prime}\big[R[Q_n;f((\cdot)_{Q^{\prime}},(\cdot)_R)]\big]
+Q^{\prime\prime}\big[R[Q_n;f((\cdot)_R,(\cdot)_{Q^{\prime\prime}})]\big].
\end{equation}
\end{thm}
See also \cite[Theorems 1 and 3]{GN2009}.

In view of the integral representation \eqref{e:2.9}, we may have
definite blending cubature formulae.

\begin{dfn}
The blending cubature formula $S_n$ is called positive (negative)
definite of order $(r,s)$ if $K_{r,s}(t,\tau)\geq 0$ (resp.
$K_{r,s}(t,\tau)\leq 0$) in $[a,b]^2$.
\end{dfn}
Definite cubature formulae of order $(r,s)$ provide one-sided
approximation to $I[f]$ for integrands from the class
$$
\mathcal{C}^{r,s}[a,b]=\{f\in C^{r,s}([a,b]^2)\,:\,D^{r,s}f\
\text{does not change sign in}\ (a,b)^2\}\,.
$$
Set
\begin{eqnarray*}
&&\mathcal{C}^{r,s}_{+}[a,b]=\{f\in \mathcal{C}^{r,s}[a,b]\,:\,
D^{r,s}f\geq 0\ \text{ in }\ [a,b]^2\},\\
&&\mathcal{C}^{r,s}_{-}[a,b]=\{f\in \mathcal{C}^{r,s}[a,b]\,:\,
-f\in \mathcal{C}^{r,s}_{+}[a,b]\}.
\end{eqnarray*}
If $S_n$ is, say, positive definite of order $(r,s)$, then
$R[S_n;f]\geq 0$ for $f\in\mathcal{C}^{r,s}_{+}[a,b]$ and the
reversed inequality holds when $f\in\mathcal{C}^{r,s}_{-}[a,b]$.

Various sufficient conditions for definiteness of blending cubature
formulae are given in \cite[Theorem~5]{GN2009}.

We conclude this section with a general observation about definite
blending cubature formulae.
\begin{thm}\label{t:4}
Let $\,(S^{\prime}\,,\,S^{\prime\prime})\,$ be a pair of positive
(negative) definite blending cubature formulae of order $(r,s)$.
Assume that, for some $c>0$, the cubature formula
$$
\widehat{S}:=(c+1)\,S^{\prime}-c\,S^{\prime\prime}
$$
is negative (positive) definite of order $(r,s)$. Then the following
inequalities hold true whenever $\,f\,$ belongs to
$\mathcal{C}^{r,s}[a,b]$:
\begin{enumerate}[~~(i)~~]
\item
$\ds{|R[S^{\prime};f]|\leq\frac{c}{c+1}\,|R[S^{\prime\prime};f]|}$\,;
\item $\ds{|R[S^{\prime};f]|\leq
c\,|S^{\prime}[f]-S^{\prime\prime}[f]|}$\,;
\item $\ds{|R[S^{\prime\prime};f]|\leq
(c+1)\,|S^{\prime}[f]-S^{\prime\prime}[f]|}$\,.
\end{enumerate}
\end{thm}
\begin{pf}
Let us consider, e.g., the case when $S^{\prime}$ and
$S^{\prime\prime}$ are negative definite and $\,\widehat{S}\,$ is
positive definite, of order $(r,s)$. Without loss of generality we
may assume that $f\in \mathcal{C}^{r,s}_{+}[a,b]$. Then
$\,R[S^{\prime};f]\leq 0$, $\,R[S^{\prime\prime};f]\leq 0$, and
$\,R[\widehat{S};f]\geq 0$, therefore
$$
0\leq
R[\widehat{S};f]=(c+1)\,R[S^{\prime};f]-c\,R[S^{\prime\prime};f]\,,
$$
and hence
$$
-R[S^{\prime};f]\leq-\frac{c}{c+1}\,R[S^{\prime\prime};f]\,,
$$
which, in this case, is claim (i) of Theorem \ref{t:4}. Claim (iii)
follows from
\[
\begin{split}
|S^{\prime}[f]-S^{\prime\prime}[f]|&=|R[S^{\prime\prime};f]-R[S^{\prime};f]|
\geq |R[S^{\prime\prime};f]|-|R[S^{\prime};f]|\\
&\geq |R[S^{\prime\prime};f]|-\frac{c}{c+1}\,|R[S^{\prime\prime};f]|
=\frac{1}{c+1}\,|R[S^{\prime\prime};f]|\,.
\end{split}
\]
Part (ii) is a consequence of (iii) and (i). The proof of the case
when $S^{\prime}$ and $S^{\prime\prime}$ are positive definite and
$\,\widehat{S}\,$ is negative definite of order $(r,s)$ is
analogous, and we omit it.\qed
\end{pf}

See \cite[Theorem 1]{AN2017} for analogous statement concerning
definite quadrature formulae.

\section{Proofs}
\subsection{Proof of Theorem~\ref{t:1.1}}
The cubature formula $S_n^{-}$ is obtained through the scheme
\eqref{e:2.7}, where $Bf$ is the Hermite blending interpolant for
$f$ at the blending grid $G(X,Y)$ formed by
$X=Y=\Big(\frac{a+b}{2},\frac{a+b}{2}\Big)$, and
$C_n[f]=C_n^{Tr}[f]$. In this case
$Q^{\prime}=Q^{\prime\prime}=Q^{Mi}$ is the midpoint quadrature
formula
$$
Q^{Mi}[g]=(b-a)g\big(\frac{a+b}{2}\big),
$$
and \eqref{e:2.13} implies the representation \eqref{e:1.1} of
$S_n^{-}$. According to \eqref{e:2.12}, the Peano kernel
$K_{2,2}(S_n^{-};t,\tau)$ is equal to
\begin{equation}\label{e:3.1}
K_2(Q^{Mi};t)K_2(Q_n^{Tr};\tau) +K_2(Q^{Mi};\tau)K_2(Q_n^{Tr};t)
-K_2(Q_n^{Tr};t)K_2(Q_n^{Tr};\tau).
\end{equation}
Since $Q_n^{Tr}$ and $Q^{Mi}$ are respectively negative definite and
positive definite quadrature formulae of order two, \eqref{e:3.1}
implies $K_{2,2}(S_n^{-};t,\tau)\leq 0$ in $[a,b]^2$, i.e.,
$S_n^{-}$ is a negative definite cubature formula of order $(2,2)$.
By the mean value theorem, there is a point $P\in [a,b]^2$ such that
$$
R[S_n^{-};f]
=D^{2,2}f(P)\iint\limits_{[a,b]^2}K_{2,2}(S_n^{-};t,\tau)\,dtd\tau.
$$
From \eqref{e:3.1} and
$$
\int\limits_{a}^{b}K_2(Q^{Mi};t)\,dt=\frac{(b-a)^{3}}{24}, \qquad
\int\limits_{a}^{b}K_2(Q_{n}^{Tr};t)\,dt=-\frac{(b-a)^{3}}{12n^2},
$$
we find
$$
\iint\limits_{[a,b]^2}\!K_{2,2}(S_n^{-};t,\tau)\,dtd\tau=
-\frac{(b-a)^{6}}{144n^2}\Big(1+\frac{1}{n^2}\Big),
$$
thus proving part (i) of Theorem~A.

Theorem~\ref{t:1.1} follows from Theorem~\ref{t:4}, applied with
$r=s=2$, $S^{\prime}=S_{2n}^{-}$, $S^{\prime\prime}=S_{n}^{-}$, and
$c=1$. We shall prove that the cubature formula
$$
\widehat{S}=(c+1)S_{2n}^{-}-c\,S_{n}^{-}
$$
is positive definite of order $(2,2)$ for every $c\geq 1$, or,
equivalently, that
\begin{equation}\label{e:3.2}
\varphi(t,\tau):=(c+1)K_{2,2}(S_{2n}^{-};t,\tau)
-c\,K_{2,2}(S_{n}^{-};t,\tau)\geq 0,\quad (t,\tau)\in [a,b]^2
\end{equation}
when $c\geq 1$. Moreover, we shall show that $c=1$ is the smallest
constant with this property. We assume in what follows that $c\geq
1$, and set
\begin{equation}\label{e:3.3}
t_{i}=\tau_{i}=a+i\,h,\quad i=0,1,\ldots,2n,\quad h=\frac{b-a}{2n}.
\end{equation}
We split $[a,b]^2$ in the following way:
$$
[a,b]^2=\bigcup_{k=0}^{n-1}\bigcup_{\ell=0}^{n-1}\Delta_{k,\ell},\quad
\Delta_{k,\ell}=\{(t,\tau)\,:\, t_{2k}\leq t\leq
t_{2k+2},\,\tau_{2\ell}\leq \tau\leq\tau_{2\ell+2}\},
$$
with a further decomposition $\ds
\Delta_{k,\ell}=\bigcup_{i=1}^{4}\Delta_{k,\ell}^{i}$, where
\begin{eqnarray*}
&&\Delta_{k,\ell}^{1}=\{(t,\tau)\,:\, t_{2k}\leq t\leq
t_{2k+1},\,\tau_{2\ell}\leq \tau\leq\tau_{2\ell+1}\},\\
&&\Delta_{k,\ell}^{2}=\{(t,\tau)\,:\, t_{2k+1}\leq t\leq
t_{2k+2},\,\tau_{2\ell}\leq \tau\leq\tau_{2\ell+1}\},\\
&&\Delta_{k,\ell}^{3}=\{(t,\tau)\,:\, t_{2k}\leq t\leq
t_{2k+1},\,\tau_{2\ell+1}\leq \tau\leq\tau_{2\ell+2}\},\\
&&\Delta_{k,\ell}^{4}=\{(t,\tau)\,:\, t_{2k+1}\leq t\leq
t_{2k+2},\,\tau_{2\ell+1}\leq \tau\leq\tau_{2\ell+2}\}.
\end{eqnarray*}
The Peano kernels $K_2(Q^{Mi};\cdot)$, $K_2(Q_n^{Tr};\cdot)$ and
$K_2(Q_{2n}^{Tr};\cdot)$, which according to \eqref{e:3.1} and
\eqref{e:3.2} occur in $\varphi$, admit the representations (see
\eqref{e:2.2}):
\begin{equation}\label{e:3.4}
\begin{split}
&K_2(Q^{Mi};t)=\frac{(t-a)^2}{2}-(b-a)\Big(t-\frac{a+b}{2}\Big)_{+},
\qquad t\in [a,b],\\
&K_2(Q_n^{Tr};t)=\frac{1}{2}\big(t-t_{2k})(t-t_{2k+2}),\quad t\in
\big[t_{2k},t_{2k+2}\big],\ \ 0\leq k\leq n-1,\\
&K_2(Q_{2n}^{Tr};t)=\frac{1}{2}\big(t-t_{i})(t-t_{i+1}),\quad t\in
\big[t_{i},t_{i+1}\big],\ \ 0\leq i\leq 2n-1.
\end{split}
\end{equation}
Since either of these Peano kernels is an even function with respect
$\frac{a+b}{2}$, i.e., $K_2(\cdot;\frac{a+b}{2}-t)=
K_2(\cdot;\frac{a+b}{2}+t)$ for $a\leq t\leq \frac{a+b}{2}$, then
$\varphi(t,\tau)$ possesses the same property in each of variables
$t$ and $\tau$. Therefore it suffices to verify \eqref{e:3.4} only
for $a\leq t,\tau\leq \frac{a+b}{2}$.

Assume first that $(t,\tau)\in\Delta_{0,0}^{1}$, then using
\eqref{e:3.1} and \eqref{e:3.4} we find {\small
\begin{equation*}
\begin{split}
\varphi(t,\tau)=\frac{(t-a)(\tau-a)}{4}\Big[&(c+1)\big[(t-a)(\tau-t_1)
+(\tau-a)(t-t_1)-(t-t_1)(\tau-t_1)\big]\\
&-c\big[(t-a)(\tau-t_2)
+(\tau-a)(t-t_2)-(t-t_2)(\tau-t_2)\big]\Big],
\end{split}
\end{equation*}}
and substituting $t_2=t_1+h$, get
$$
\varphi(t,\tau)=\frac{(t-a)(\tau-a)}{4}\Big[(t-a)(\tau-t_1)
+(\tau-a)(t-t_1)-(t-t_1)(\tau-t_1)+3c\,h^2\Big].
$$
Since the expression in the brackets is a bilinear function, its
minimal value in $\Delta_{0,0}^{1}$ is attained at a vertex of
$\Delta_{0,0}^{1}$. Hence, thus this expression is non-negative when
$c\geq 1/3$, and consequently $\varphi(t,\tau)\geq 0$ for
$(t,\tau)\in\Delta_{0,0}^{1}$.

Next, we compare the values of $\varphi(t,\tau)$ in
$\Delta_{k,\ell}^{1}$ and $\Delta_{k+1,\ell}^{1}$, assuming that
$n\geq 2$ and $\ell\leq\frac{n-1}{2}$, $k\leq\frac{n-2}{2}$. By
using \eqref{e:3.4}, for $(t,\tau)\in\Delta_{k,\ell}^{1}$ we obtain
\begin{equation*}
\begin{split}
4\varphi(t,\tau)=&(c+1)\Big[(t-a)^2(\tau-t_{2\ell})(\tau-t_{2\ell+1})
+(\tau-a)^2(t-t_{2k})(t-t_{2k+1})\\
&\qquad\qquad -(t-t_{2k})(t-t_{2k+1})(\tau-t_{2\ell})(\tau-t_{2\ell+1})\Big]\\
&-c\Big[(t-a)^2(\tau-t_{2\ell})(\tau-t_{2\ell+2})
+(\tau-a)^2(t-t_{2k})(t-t_{2k+2})\\
&\qquad\qquad -
(t-t_{2k})(t-t_{2k+2})(\tau-t_{2\ell})(\tau-t_{2\ell+2})\Big].
\end{split}
\end{equation*}
With $t=t_{2k}+u\,h$, $\tau=t_{2\ell}+v\,h$, $\;(u,v)\in [0,1]^2$,
this expression simplifies to
\begin{equation}\label{e:3.5}
\begin{split}
\psi_{k,\ell}(u,v)=h^4\Big[&(2k+u)^2 v(v-1+c)+(2\ell+v)^2u(u-1+c)\\
&-uv(1-u)(1-v)+c\,uv(3-u-v)\big)\Big].
\end{split}
\end{equation}
Analogously, by substituting $t=t_{2k+2}+u\,h$,
$\tau=t_{2\ell}+v\,h$ with $(u,v)\in [0,1]^2$ we find that
$4\varphi(t,\tau)=\psi_{k+1,\ell}(u,v)$ for $(t,\tau)\in
\Delta_{k+1,\ell}^{1}$. Since
$$
\psi_{k+1,\ell}(u,v)-\psi_{k,\ell}(u,v)=4h^4(2k+1+u)v(c-1+v)\geq
0,\qquad (u,v)\in [0,1]^2,
$$
the following implication holds true for $0\leq k\leq\frac{n-2}{2}$
and $0\leq\ell\leq\frac{n-1}{2}$:
\begin{equation}\label{e:3.6}
\varphi(t,\tau)\geq 0,\ (t,\tau)\in\Delta_{k,\ell}^{1}\
\Longrightarrow\ \varphi(t,\tau)\geq 0,\
(t,\tau)\in\Delta_{k+1,\ell}^{1}.
\end{equation}

Next, we compare the values of $\varphi$ on $\Delta_{k,\ell}^{1}$
and $\Delta_{k,\ell}^{2}$, where $0\leq k,\ell\leq\frac{n-1}{2}$.
For $(t,\tau)\in\Delta_{k,\ell}^{2}$ we have
\begin{equation*}
\begin{split}
4\varphi(t,\tau)=&(c+1)\Big[(t-a)^2(\tau-t_{2\ell})(\tau-t_{2\ell+1})
+(\tau-a)^2(t-t_{2k+1})(t-t_{2k+2})\\
&\qquad\qquad -(t-t_{2k+1})(t-t_{2k+2})(\tau-t_{2\ell})(\tau-t_{2\ell+1})\Big]\\
&-c\Big[(t-a)^2(\tau-t_{2\ell})(\tau-t_{2\ell+2})
+(\tau-a)^2(t-t_{2k})(t-t_{2k+2})\\
&\qquad\qquad -
(t-t_{2k})(t-t_{2k+2})(\tau-t_{2\ell})(\tau-t_{2\ell+2})\Big].
\end{split}
\end{equation*}
With $t=t_{2k+2}-u\,h$, $\tau=t_{2\ell}+v\,h$, where $(u,v)\in
[0,1]^2$, the above expression becomes
\begin{equation*}
\begin{split}
\widetilde{\psi}_{k,\ell}(u,v)=h^4\Big[&(2k+2-u)^2 v(v-1+c)+(2\ell+v)^2u(u-1+c)\\
&-uv(1-u)(1-v)+c\,uv(3-u-v)\big)\Big].
\end{split}
\end{equation*}
Since
$\widetilde{\psi}_{k,\ell}(u,v)-\psi_{k,\ell}(u,v)=4h^4(2k+1-u)(1-u)v(v-1+c)\geq
0$ for $(u,v)\in [0,1]^2$, we have the following implication when
$0\leq k,\ell\leq\frac{n-1}{2}$:
\begin{equation}\label{e:3.7}
\varphi(t,\tau)\geq 0,\ (t,\tau)\in\Delta_{k,\ell}^{1}\
\Longrightarrow\ \varphi(t,\tau)\geq 0,\
(t,\tau)\in\Delta_{k,\ell}^{2}.
\end{equation}
We apply the same procedure to deduce similar conclusions in the
"vertical" direction. Thus, for $0\leq k,\ell\leq\frac{n-1}{2}$ we
have the implication
\begin{equation}\label{e:3.8}
\varphi(t,\tau)\geq 0,\ (t,\tau)\in\Delta_{k,\ell}^{1}\
\Longrightarrow\ \varphi(t,\tau)\geq 0,\
(t,\tau)\in\Delta_{k,\ell}^{3},
\end{equation}
and for $k\leq\frac{n-1}{2}$, $\ell\leq\frac{n-2}{2}$,
\begin{equation}\label{e:3.9}
\varphi(t,\tau)\geq 0,\ (t,\tau)\in\Delta_{k,\ell}^{1}\
\Longrightarrow\ \varphi(t,\tau)\geq 0,\
(t,\tau)\in\Delta_{k,\ell+1}^{1}.
\end{equation}
Since we already showed that $\varphi(t,\tau)\geq 0$ in
$\Delta_{0,0}^{1}$, from \eqref{e:3.6}, \eqref{e:3.7}, \eqref{e:3.8}
and \eqref{e:3.9} we conclude that $\varphi(t,\tau)\geq 0$ for all
$(t,\tau)$ satisfying $a\leq t,\,\tau\leq\frac{a+b}{2}$. By
symmetry, it follows that $\varphi(t,\tau)\geq 0$ in $[a,b]^2$.

Thus, we showed that the cubature formula
$$
\widehat{S}=(c+1)S_{2n}^{-}-c\,S_{n}^{-}
$$
is positive definite of order $(2,2)$ whenever $c\geq 1$.
Theorem~\ref{t:1.1} follows from Theorem~\ref{t:4}, applied with
$r=s=2$, $S^{\prime}=S_{2n}^{-}$, $S^{\prime\prime}=S_{n}^{-}$, and
$c=1$. \qed

\begin{rmk}
To demonstrate that $\widehat{S}$ is not positive definite when
$c<1$, we show that $K_{2,2}(\widehat{S};t,\tau)$ assumes negative
values on the segment $t=\tau$ in $\Delta_{k,k}^{1}$, $1\leq k\leq
\frac{n-1}{2}$. In view of \eqref{e:3.5}, the trace of
$K_{2,2}(\widehat{S};\cdot)$ on this segment is the univariate
function
$$
g(u)=h^4\big[2(2k+u)^2u(u-1+c)-u^2(1-u)^2+c\,u^2(3-2u)\big],\qquad
u\in [0,1].
$$
From $g(0)=0$ and $g^{\prime}(0)=8(c-1)h^4k^2$ it follows that
$g(u)$ assumes negative values in $(0,1)$ when $c<1$ .
\end{rmk}

\subsection{Proof of Theorem~\ref{t:1.2}}
$S_n^{+}$ is built according to \eqref{e:2.7}, where $Bf$ is the
Lagrange blending interpolant for $f$ at the blending grid $G(X,Y)$
with $X=Y=\{a,b\}$, and $C_n[f]=C_n^{Tr}[f]$. In this case
$Q^{\prime}=Q^{\prime\prime}=Q^{Tr}$, and the representation
\eqref{e:1.2} of $S_n^{+}$ follows from \eqref{e:2.13}.

From \eqref{e:2.11} we have
\begin{equation}\label{e:3.10}
K_{2,2}(S_n^{+};t,\tau)=K_2(Q^{Tr};\tau)K_2(Q_n^{Tr};t)
+K_2(Q_n^{Tr};\tau)Q_n^{Tr}[\mathcal{K}_2(\cdot,t)],
\end{equation}
where
$$
\mathcal{K}_2(x,t)=(x-t)_{+}-\frac{x-a}{b-a}\,(b-t).
$$
Since $Q_n^{Tr}$ (in particular, $Q^{Tr}$) is negative definite of
order two and $\mathcal{K}_2(x,t)\leq 0$ for $t,\tau\in [a,b]$, it
follows from \eqref{e:3.10} that $K_{2,2}(S_n^{+};t,\tau)\geq 0$ in
$[a,b]^2$, i.e., $S_n^{+}$ is a positive definite cubature formula
of order $(2,2)$. According to \eqref{e:2.12},
\begin{equation}\label{e:3.11}
\begin{split}
K_{2,2}(S_n^{+};t,\tau)=&K_2(Q^{Tr};t)K_2(Q_n^{Tr};\tau)
+K_2(Q^{Tr}\!;\tau)K_2(Q_n^{Tr};t)\\
&-K_2(Q_n^{Tr};t)K_2(Q_n^{Tr};\tau).
\end{split}
\end{equation}
Claim (ii) of Theorem~A follows from the mean value theorem, applied
to the integral representation of the remainder $R[S_n^{+};f]$,
taking into account \eqref{e:3.11}.

The proof of Theorem~\ref{t:1.2} is based on Theorem~\ref{t:4}
applied with $r=s=2$, $S^{\prime}=S_{2n}^{+}$ and
$S^{\prime\prime}=S_{n}^{+}$. We shall find condition on the
constant $c>0$ so that the cubature formula
$\widehat{S}=(c+1)\,S_{2n}^{+}-c\,S_{2n}^{+}$ is negative definite
of order $(2,2)$, i.e.,
\begin{equation}\label{e:3.12}
\varphi(t,\tau):=(c+1)K_{2,2}(S_{2n}^{+};t,\tau)
-c\,K_{2,2}(S_{n}^{+};t,\tau)\leq 0,\quad (t,\tau)\in [a,b]^2
\end{equation}
The verification of \eqref{e:3.12} goes along the same lines as that
of \eqref{e:3.2}: with $\{t_i\}_{i=0}^{2n}$ defined in \eqref{e:3.3}
we decompose $[a,b]^2$ in the same way, and study $\varphi(t,\tau)$
on the different sub-domains. By symmetry, we may restrict our
considerations to the region $a\leq t,\tau\leq\frac{a+b}{2}$. Assume
first that $(t,\tau)\in\Delta_{0,0}^{1}$. Using \eqref{e:3.11},
\eqref{e:3.4} and
$$
K_2(Q^{Tr};t)=\frac{1}{2}(t-a)(t-b),
$$
we find that for $(t,\tau)\in \Delta_{0,0}^1$ the function
$\varphi(t,\tau)$ takes the form {\small
\begin{equation*}
\begin{split}
\varphi(t,\tau)=\frac{(t-a)(\tau-a)}{4}\Big[&(c+1)\big[(t-b)(\tau-t_1)
+(\tau-b)(t-t_1)-(t-t_1)(\tau-t_1)\big]\\
&-c\big[(t-b)(\tau-t_2)
+(\tau-b)(t-t_2)-(t-t_2)(\tau-t_2)\big]\Big].
\end{split}
\end{equation*}}
Replacing $t_2=t_1+h$ and using that $b-a=2n\,h$, after
simplification we get
$$
\varphi(t,\tau)=\frac{(t-a)(\tau-a)}{4}\Big[t\tau-b(t+\tau)+t_1(2b-t_1)
-c(4n-3)h^2\Big].
$$
The expression in the brackets is a bilinear function in $t$ and
$\tau$ and attains its extreme values at the vertices of
$\Delta_{0,0}^1$. Its maximum in $\Delta_{0,0}^1$ is attained at
$(t,\tau)=(t_0,t_0)=(a,a)$ and equals $\big(4n-1-(4n-3)c\big)h^2$.
Hence, the inequality
\begin{equation}\label{e:3.13}
c\geq\frac{4n-1}{4n-3}
\end{equation}
is a necessary and sufficient condition for $\varphi(t,\tau)\leq 0$
in $\Delta_{0,0}^{1}$.

The rest of the proof of \eqref{e:3.12} follows the scheme of the
proof of \eqref{e:3.2}: under the assumption \eqref{e:3.13} we
establish the implications
\begin{equation}\label{e:3.14}
\varphi(t,\tau)\leq 0,\ (t,\tau)\in\Delta_{k,\ell}^{1}\
\Longrightarrow\ \varphi(t,\tau)\leq 0,\
(t,\tau)\in\Delta_{k+1,\ell}^{1}
\end{equation}
for $0\leq k\leq\frac{n-2}{2}$ and $0\leq \ell\leq\frac{n-1}{2}$;
\begin{equation}\label{e:3.15}
\varphi(t,\tau)\leq 0,\ (t,\tau)\in\Delta_{k,\ell}^{1}\
\Longrightarrow\ \varphi(t,\tau)\leq 0,\
(t,\tau)\in\Delta_{k,\ell}^{2}
\end{equation}
for $0\leq k,\ell\leq\frac{n-1}{2}$;
\begin{equation}\label{e:3.16}
\varphi(t,\tau)\leq 0,\ (t,\tau)\in\Delta_{k,\ell}^{1}\
\Longrightarrow\ \varphi(t,\tau)\leq 0,\
(t,\tau)\in\Delta_{k,\ell+1}^{1}
\end{equation}
for $0\leq k\leq\frac{n-1}{2}$ and $0\leq \ell\leq\frac{n-2}{2}$,
and
\begin{equation}\label{e:3.17}
\varphi(t,\tau)\leq 0,\ (t,\tau)\in\Delta_{k,\ell}^{1}\
\Longrightarrow\ \varphi(t,\tau)\leq 0,\
(t,\tau)\in\Delta_{k,\ell}^{3}
\end{equation}
for $0\leq k,\ell\leq\frac{n-1}{2}$.

Let us prove \eqref{e:3.14}. For $(t,\tau)\in\Delta_{k,\ell}^{1}$ we
have {\small
\begin{equation*}
\begin{split}
4\varphi(t,\tau)=&(c+1)\Big[(t-a)(t-b)(\tau-t_{2\ell})(\tau-t_{2\ell+1})
+(\tau-a)(\tau-b)(t-t_{2k})(t-t_{2k+1})\\
&\qquad\qquad -(t-t_{2k})(t-t_{2k+1})(\tau-t_{2\ell})(\tau-t_{2\ell+1})\Big]\\
&-c\Big[(t-a))t-b)(\tau-t_{2\ell})(\tau-t_{2\ell+2})
+(\tau-a)(\tau-b)(t-t_{2k})(t-t_{2k+2})\\
&\qquad\qquad -
(t-t_{2k})(t-t_{2k+2})(\tau-t_{2\ell})(\tau-t_{2\ell+2})\Big].
\end{split}
\end{equation*}}
With $t=t_{2k}+u\,h$, $\tau=t_{2\ell}+v\,h$, $\;(u,v)\in [0,1]^2$,
we have $4\varphi(t,\tau)=h^4\psi_{k,\ell}(u,v)$, where {\small
\begin{equation*}
\begin{split}
\psi_{k,\ell}(u,v)=&(2k+u)(2k+u-2n) v(v-1+c)+(2\ell+v)(2\ell+v-2n)u(u-1+c)\\
&-uv(1-u)(1-v)+c\,uv(3-u-v)\big).
\end{split}
\end{equation*}}
For $(t,\tau)\in\Delta_{k+1,\ell}^{1}$, we substitute
$t=t_{2k+2}+u\,h$, $\tau=t_{2\ell}+v\,h$, $\;(u,v)\in [0,1]^2$, and
obtain $4\varphi(t,\tau)=\psi_{k+1,\ell}(u,v)$. Now \eqref{e:3.14}
follows from
$$
\psi_{k+1,\ell}(u,v)-\psi_{k,\ell}(u,v)=4v(v-1+c)(2k+1+u-n)\leq 0.
$$
The proof of implications \eqref{e:3.15}, \eqref{e:3.16} and
\eqref{e:3.17} is similar, therefore we skip it. Since
$\varphi(t,\tau)\leq 0$ in $\Delta_{0,0}^{1}$ for every $c$
satisfying \eqref{e:3.13}, it follows from \eqref{e:3.14},
\eqref{e:3.15}, \eqref{e:3.16} and \eqref{e:3.17} that
$\varphi(t,\tau)\leq 0$ whenever $a\leq t,\tau\leq\frac{a+b}{2}$. By
symmetry, $\varphi(t,\tau)\leq 0$ for every $(t,\tau)\in [a,b]^2$.
Hence, the cubature formula
$$
\widehat{S}=(c+1)\,S_{2n}^{+}-c\,S_{2n}^{+}
$$
is negative definite of order $(2,2)$ whenever $c$ satisfies
\eqref{e:3.13}. Theorem~\ref{t:1.2} now follows from
Theorem~\ref{t:4} applied with $r=s=2$, $S^{\prime}=S_{2n}^{+}$,
$S^{\prime\prime}=S_{n}^{+}$ and $c=\frac{4n-1}{4n-3}$. Let us point
out that this value of $c$ is the best possible, i.e., the smallest
one for which $\widehat{S}$ is negative definite cubature formula of
order $(2,2)$.
\section{Numerical example}
We illustrate theoretical results provided by Theorems~\ref{t:1.1}
and \ref{t:1.2} by considering double integrals over $[0,1]^2$ of
the functions $f(x,y)=e^{xy}$ and $g(x,y)=\sin(xy)$. The numerical
values $I[f]\approx 1.317902151$ and $I[g]\approx 0.239811742$ are
found from
$$
\iint\limits_{[0,1]^2}e^{xy}\,dxdy=\int\limits_{0}^{1}\frac{e^x-1}{x}\,dx,
\qquad
\iint\limits_{[0,1]^2}\sin({xy})\,dxdy=\int\limits_{0}^{1}\frac{1-\cos
x}{x}\,dx
$$
and appropriate truncation of the power series of the integrands in
the univariate integrals.
\begin{table}[h]\label{tab:1}
\begin{center}
{\scriptsize
\begin{tabular}{| c | c|  c | c|  c|}
\hline
& & & &\\
n & $\displaystyle R[S_n^-;f]$ &
$\displaystyle\frac{1}{2}\big|[S_{2n}^-[f]-S_n^-[f]\big|$ &
$\displaystyle R[S_n^+;f]$ &
$\displaystyle\frac{4n-1}{4n-3}\big|S_{2n}^+[f]-S_n^+[f]\big|$
\\
& & & &\\
\hline
& & & &\\
 4 & $-1.947\times 10^{-3}$ &  & $3.615\times 10^{-3}$ &\\
  &  & $7.411\times 10^{-4}$ &  & $3.101\times 10^{-3}$\\
 8 & $-4.648\times 10^{-4}$ &  & $9.274\times 10^{-4}$&\\
  &  & $1.750\times 10^{-4}$ &  & $7.419\times 10^{-4}$\\
 16 & $-1.148\times 10^{-4}$ &  & $2.333\times 10^{-4}$&\\
  &  & $4.310\times 10^{-5}$ &  & $1.806\times 10^{-4}$\\
 32 & $-2.862\times 10^{-5}$ &  & $5.842\times 10^{-5}$&\\
  &  & $1.073\times 10^{-5}$ &  & $4.451\times 10^{-5}$\\
 64 & $-7.149\times 10^{-6}$ &  & $1.461\times 10^{-5}$&\\
  &  & $2.681\times 10^{-6}$ &  & $1.104\times 10^{-5}$\\
 128 & $-1.787\times 10^{-6}$ &  & $3.653\times 10^{-6}$&\\
& & & &\\
\hline
\end{tabular} }
\end{center}
\caption{Remainders and error bounds of $S_n^{-}[f]$ and
$S_n^{+}[f]$, $f(x,y)=e^{xy}$.}
\end{table}

\begin{table}[h]\label{tab:2}
\begin{center}
{\scriptsize
\begin{tabular}{| c | c | c | c | c |}
\hline
& & & &\\
n & $\displaystyle R[S_n^-;g]$ &
$\displaystyle\frac{1}{2}\big|[S_{2n}^-[g]-S_n^-[g]\big|$ &
$\displaystyle R[S_n^+;g]$
& $\displaystyle\frac{4n-1}{4n - 3}\big|S_{2n}^+[g]-S_n^+[g]\big|$\\
& & & &\\
\hline
& & & &\\
 4 & $6.300\times 10^{-4}$ &  & $-1.129\times 10^{-3}$ &\\
  &  & $2.397\times 10^{-4}$ &  & $9.697\times 10^{-4}$\\
 8 & $1.507\times 10^{-4}$ &  & $-2.886\times 10^{-4}$ &\\
  &  & $5.674\times 10^{-5}$ &  & $2.309\times 10^{-4}$\\
 16 & $3.726\times 10^{-5}$ &  & $-7.254\times 10^{-5}$ &\\
  &  & $1.399\times 10^{-5}$ &  & $5.616\times 10^{-5}$\\
 32 & $9.289\times 10^{-6}$ &  & $-1.816\times 10^{-5}$ &\\
  &  & $3.484\times 10^{-6}$ &  & $1.384\times 10^{-5}$\\
 64 & $2.321\times 10^{-6}$ &  & $-4.541\times 10^{-6}$ &\\
  &  & $8.703\times 10^{-7}$ &  & $3.433\times 10^{-6}$\\
 128 & $5.801\times 10^{-7}$ &  & $-1.135\times 10^{-6}$ &\\
& & & &\\
\hline
\end{tabular} }
\end{center}
\caption{Remainders and error bounds of $S_n^{-}[g]$ and
$S_n^{+}[g]$, $g(x,y)=\sin(xy)$.}
\end{table}

Table 1 and Table 2 present the remainders and a posteriori error
estimates of cubature formulae $S_n^-$ and $S_n^+$ applied to
$f(x,y)$ and $g(x,y)$, respectively. The remainder signs
($R[S_n^-;f]<0$, $R[S_n^+;f]>0$, $R[S_n^-;g]>0$, $R[S_n^+;g]<0$) are
in agreement with the fact that $D^{2,2}f>0$ and $D^{2,2}g<0$ in
$[0,1]^2$. In both tables it is observed  that a posteriori error
estimates exceed the true error magnitude by a factor ranging
between $1.5$ and $1.6$ for $S_n^{-}$ and between $3.02$ and $3.36$
for $S_n^{+}$. Of course, the mean cubature formula
$(S_n^{+}+S_n^{-})/2$ furnishes another reasonable approximation to
$I$, and for integrands from $\mathcal{C}^{2,2}[0,1]$ it has error
magnitude not exceeding $|R[S_n^{+};\cdot]+R[S_n^{-};\cdot]|/2$.
\begin{rmk}
One may wonder if similar results may be established for the
modified composite midpoint rules, constructed through the scheme
\eqref{e:2.7} with blending interpolation on the same blending
grids. The answer is in the affirmative regarding the definiteness,
i.e., the analogue of Theorem~A. Unfortunately, we cannot prove the
analogues of Theorems~\ref{t:1.1} and \ref{t:1.2} (monotonicity and
a posteriori error estimates), since Theorem~\ref{t:4} is not
applicable in this situation, due to the fact that the associated
Peano kernels $K_{2,2}(S_n;\cdot)$ and $K_{2,2}(S_{2n};\cdot)$
coincide on some subregions of $[a,b]^2$.
\end{rmk}

\end{document}